\pdfoutput=1
\RequirePackage{ifpdf}
\ifpdf 
\documentclass[pdftex]{sigma}
\else
\documentclass{sigma}
\fi

\newtheorem*{question*}{Question}

\begin{document}

\allowdisplaybreaks

\newcommand{\arXivNumber}{1407.8291}

\renewcommand{\PaperNumber}{112}

\FirstPageHeading

\ShortArticleName{Conf\/igurations of Points and the Symplectic Berry--Robbins Problem}

\ArticleName{Conf\/igurations of Points\\
and the Symplectic Berry--Robbins Problem}

\Author{Joseph MALKOUN}
\AuthorNameForHeading{J.~Malkoun}

\Address{Department of Mathematics and Statistics, Notre Dame University-Louaize, Lebanon}
\Email{\href{mailto:joseph.malkoun@ndu.edu.lb}{joseph.malkoun@ndu.edu.lb}}

\ArticleDates{Received August 23, 2014, in f\/inal form December 17, 2014; Published online December 19, 2014}

\Abstract{We present a~new problem on conf\/igurations of points, which is a~new version of a~similar problem by Atiyah
and Sutclif\/fe, except it is related to the Lie group $\operatorname{Sp}(n)$, instead of the Lie group $\operatorname{U}(n)$.
Denote by $\mathfrak{h}$ a~Cartan algebra of $\operatorname{Sp}(n)$, and~$\Delta$ the union of the zero sets of the
roots of $\operatorname{Sp}(n)$ tensored with $\mathbb{R}^3$, each being a~map from $\mathfrak{h} \otimes \mathbb{R}^3
\to \mathbb{R}^3$.
We wish to construct a~map $(\mathfrak{h} \otimes \mathbb{R}^3) \backslash \Delta \to \operatorname{Sp}(n)/T^n$ which is
equivariant under the action of the Weyl group $W_n$ of $\operatorname{Sp}(n)$ (the symplectic Berry--Robbins problem).
Here, the target space is the f\/lag manifold of $\operatorname{Sp}(n)$, and $T^n$ is the diagonal~$n$-torus.
The existence of such a~map was proved by Atiyah and Bielawski in a~more general context.
We present an explicit smooth candidate for such an equivariant map, which would be a~genuine map provided a~certain
linear independence conjecture holds.
We prove the linear independence conjecture for $n=2$.}


\Classification{51F99; 17B22}

\section{Introduction}

We introduce the relevant manifolds.
\begin{gather*}
C_n = \big\{(x_1,\dots,x_n) \in \big(\mathbb{R}^3\backslash\{0\}\big)^n;~
x_r \neq x_s~\text{and}~x_r \neq -x_s~\text{for all}~1 \leq r<s \leq n \big\}.
\end{gather*}
We remark that $C_n = (\mathfrak{h} \otimes \mathbb{R}^3) \backslash \Delta$, where $\mathfrak{h}$ is a~Cartan algebra
of $\operatorname{Sp}(n)$, and~$\Delta$ the union of the zero sets of the roots of $\operatorname{Sp}(n)$ tensored with
$\mathbb{R}^3$, each being a~map from $\mathfrak{h} \otimes \mathbb{R}^3 \to \mathbb{R}^3$.
We denote by $F_n$ the following f\/lag manifold
\begin{gather*}
F_n = \operatorname{Sp}(n)/T^n,
\end{gather*}
where $T^n$ is the diagonal~$n$-dimensional torus in $\operatorname{Sp}(n)$.
Let $W_n$ be the Weyl group of the Lie group $\operatorname{Sp}(n)$.
It is well known that $W_n$ is a~semidirect product
\begin{gather*}
W_n = (\mathbb{Z}/(2))^n \rtimes \Sigma_n,
\end{gather*}
where $\Sigma_n$ is the symmetric group on~$n$ elements, and where we think of $\mathbb{Z}/(2)$ as $\{-1,1\}$ with
multiplication.
The group $W_n$ acts on $C_n$ as follows.
The element $(1,\ldots,-1,\ldots,1)\in (\mathbb{Z}/(2))^n$, with a~$-1$ in the~$r$th position only, acts on
$(x_1,\ldots,x_n)$ by replacing $x_r$ with $-x_r$, and leaving all other $x_s$ invariant.
An element $\sigma \in \Sigma_n$ acts by permuting the~$n$ coordinates.

On the other hand, the action of $W_n$ on $F_n$ can be described as follows: an element in $(\mathbb{Z}/(2))^n$ having
$-1$ only in the~$r$th position, multiplies the~$r$th column of each point in $F_n$ by the quaternionic structure~$j$
(leaving the other columns invariant), while a~permutation~$\sigma$ simply permutes the columns of each point $gT^n \in
F_n$.

Having described the main players in our story, we consider the following question, which is a~special case of
a~question asked and solved by Atiyah and Bielawski for an arbitrary Lie group in~\cite{Atiyah-Bielawski2002}.
Another special case, for the unitary groups, was considered earlier by Berry and Robbins in~\cite{BR1997}.

\begin{question*}
Is there for each $n \geq 2$, a~continuous map $f_n: C_n \to F_n$ which is equivariant for the action of $W_n$?
\end{question*}

Actually, as we wrote earlier, Atiyah and Bielawski have already posed and solved in~\cite{Atiyah-Bielawski2002} a~more
general problem for any compact Lie group~$G$ (in our case, $G = \operatorname{Sp}(n)$).
But their solution is non-elementary, as it relies on an analysis of the Nahm equations.
Here we propose, similar to Atiyah~\cite{Atiyah2000,Atiyah2001}), and Atiyah and Sutclif\/fe~\cite{Atiyah-Sutcliffe2002},
a~more elementary construction in the same spirit as those papers, but for the case $G = \operatorname{Sp}(n)$, instead
of $G=\operatorname{U}(n)$.

\section{The main construction}

The stereographic projection is a~map $s: S^2 \to \mathbb{C} \cup \{\infty \}$.
Denoting by~$N$ the point with coordinates $(0,0,1) \in S^2$, we def\/ine the stereographic projection $s(x,y,z)$ of
a~point $(x,y,z)\in S^2 \backslash \{N\}$~by
\begin{gather*}
s(x,y,z) = \frac{x+iy}{1-z}
\end{gather*}
and we def\/ine $s(N) = \infty$.
We f\/irst associate to each conf\/iguration $\mathbf{x} \in C_n$,~$n$ polynomials $p_1,\ldots,p_n$ of $t\in \mathbb{C}$ of
degree less than or equal to $2n-1$, each def\/ined up to a~scalar factor only.
Namely, we let $p_r$ be a~polynomial having as roots the stereographic projections of
\begin{gather*}
-\frac{x_r}{\|x_r\|},
\qquad
\frac{-x_r+x_s}{\|-x_r+x_s\|}
\qquad
\text{and}
\qquad
\frac{-x_r-x_s}{\|-x_r-x_s\|}
\end{gather*}
for all $s \neq r$ ($1 \leq s \leq n$), where $\|-\|$ 
denotes the Euclidean norm on $\mathbb{R}^3$.
Similarly, we introduce~$n$ other polynomials $q_1,\ldots,q_n$, with $q_r$ having as roots the antipodals of the roots
of~$p_r$, namely the stereographic projections of $x_r/\|x_r\|$, $(x_r+x_s)/\|x_r+x_s\|$ and $(x_r-x_s)/\|x_r-x_s\|$,
for $s \neq r$.

The space of complex polynomials of degree less than or equal to $2n-1$ has a~natural structure of a~left vector space
over the quaternions $\mathbb{H} = \mathbb{C} \oplus j\mathbb{C}$, where~$j$ acts as follows.
Writing such a~polynomial as
\begin{gather*}
p(t) = \sum\limits_{r=0}^{2n-1} a_r t^r
\end{gather*}
the quaternionic structure~$j$ maps~$p$ to $jp$, def\/ined~by
\begin{gather*}
(jp)(t) = \sum\limits_{r=0}^{2n-1} a'_r t^r,
\end{gather*}
where
\begin{gather*}
a'_r = (-1)^{r+1} \bar{a}_{2n-1-r}.
\end{gather*}
Moreover, since the roots of $q_r$ are the antipodals of those of $p_r$, it follows that $q_r$ is a~constant times
$jp_r$.

If we think of the space of polynomials above with its quaternionic structure as $\mathbb{H}^n$, then the $p_r$
def\/ine~$n$ column vectors in $\mathbb{H}^n$.
In this way we assign to any point in $C_n$ an~$n$-tuple of elements of $\mathbb{H}^n$ def\/ined up to the diagonal action
of $(\mathbb{C}^*)^n$.

Suppose now that for every $x \in C_n$ the polynomials $p_1,\ldots,p_n$ are $\mathbb{H}$-linearly independent.
Then we get a~map from $C_n$ into $\operatorname{GL}(n,\mathbb{H})/(\mathbb{C}^*)^n$.
Note that $W_n$ acts on $\operatorname{GL}(n,\mathbb{H})/(\mathbb{C}^*)^n$ as follows.
An element $(1,\ldots,-1,\ldots,1)$ in the $(\mathbb{Z}/(2))^n$ subgroup of $W_n$, with $-1$ only in the~$r$th position,
has the ef\/fect of mapping the equivalence class of a~matrix $A \in \operatorname{GL}(n,\mathbb{H})$ to the equivalence
class of the matrix obtained from~$A$ by replacing its~$r$th column by its image under~$j$.
On the other hand, an element~$\sigma$ in the $\Sigma_n$ subgroup of $W_n$ maps the equivalence class of~$A$ to the
equivalence class of the matrix obtained from~$A$ by permuting its columns via~$\sigma$.
With this action, the map becomes $W_n$-equivariant, since replacing $x_r$ by $-x_r$ has the ef\/fect of replacing $p_r$
by $jp_r$, and a~permutation of the~$n$ coordinates of~$x$ corresponds to the same permutation of the~$n$ polynomials
$p_1,\ldots,p_n$.
If we follow this map by a~$W_n$-equivariant map from $\operatorname{GL}(n,\mathbb{H})/(\mathbb{C}^*)^n$ into
$\operatorname{Sp}(n)/T^n$, then this solves the symplectic Berry--Robbins problem (assuming that $p_1,\ldots, p_n$ are
$\mathbb{H}$-linearly independent).
Such a~map can be obtained for example using the quaternionic polar decomposition.
More specif\/ically, the map from $\operatorname{GL}(n,\mathbb{H})$ into $\operatorname{Sp}(n)$ mapping~$A$ to
$A(\sqrt{A^*A})^{-1}$, where $*$ denotes the quaternionic conjugate transpose, descends to a~$W_n$-equivariant map from
$\operatorname{GL}(n,\mathbb{H})/(\mathbb{C}^*)^n$ into $\operatorname{Sp}(n)/T^n$.

Motivated by the above discussion, we make the following conjecture.

\begin{conjecture}
Given any $\mathbf{x} \in C_n$, the $2n$ polynomials $p_r$ and $q_r$ $(1 \leq r \leq n)$ are linearly independent over
$\mathbb{C}$.
\end{conjecture}

Note that the conjecture is equivalent to the $\mathbb{H}$-linear independence of $p_1,\ldots,p_n$.
We remark further that our construction is also naturally equivariant under the action of $\operatorname{SO}(3)$, which
acts on $\mathbb{R}^3$ in the natural way (which is equivalent to the adjoint action of $\operatorname{SO}(3)$ on its
Lie algebra~$\mathfrak{so}(3)$), and acts on the polynomials $p_r$ and $q_r$, which are only def\/ined up to a~(non-zero)
scalar factor, via the induced action of its double cover $\operatorname{SU}(2)$ on $S^{2n-1}H$, where~$H$ is the
natural irreducible complex two-dimensional representation space of $\operatorname{SU}(2)$.
In the following section, we def\/ine a~natural determinant function, and then, we prove the conjecture above for $n=2$.

\section{A determinant function}

Let $t^{\pm}_r$ and $t^{\pm \pm}_{rs} \in \mathbb{C}P^1$ be the stereographic projections of the normalizations of the
following vectors, respectively
\begin{gather*}
\pm x_r \quad(1\leq r \leq n) \qquad\text{and}\qquad \pm x_r \pm x_s \quad(1 \leq r,s \leq n\quad\text{and}\quad r \neq s).
\end{gather*}
The Hopf map $h: \mathbb{C}^2 \backslash \{0\} \to \mathbb{C}P^1$ is the natural projection map mapping $v \in
\mathbb{C}^2 \backslash \{0\}$ to its $\mathbb{C}^*$-orbit, where $\mathbb{C}^*$ acts on $\mathbb{C}^2 \backslash \{0\}$
by scalar multiplication.

We then choose lifts $\mathbf{u}^{\pm}_r=(u^{\pm}_r, v^{\pm}_r)$ and $\mathbf{u}^{\pm \pm}_{rs}=(u^{\pm \pm}_{rs},
v^{\pm \pm}_{rs}) \in \mathbb{C}^2 \backslash \{0\}$ of these roots under the Hopf map~$h$.
Since $\mathbf{u}^{-\sigma,-\tau}_{sr}$ and $\mathbf{u}^{\sigma,\tau}_{rs}$ are both lifts of the same point, we require
that they be equal.
We then form the polynomials $p_r$ and $q_r$ as the following products
\begin{gather*}
p_r(t)  = (u^{-}_r t - v^{-}_r) \prod\limits_{s\neq r}(u^{-+}_{rs}t -v^{-+}_{rs}) \prod\limits_{s\neq r}(u^{--}_{rs} t-v^{--}_{rs}),
\\
q_r(t)  = (u^{+}_r t - v^{+}_r) \prod\limits_{s\neq r}(u^{++}_{rs}t -v^{++}_{rs}) \prod\limits_{s\neq r}(u^{+-}_{rs} t-v^{+-}_{rs}).
\end{gather*}
Thus in particular, once the lifts are chosen, the polynomials $p_r$ and $q_r$ are determined uniquely, in other words, the scalar factor of each of them gets f\/ixed.
We now form the complex $2n$ by $2n$ matrix
\begin{gather*}
M = (p_1, q_1, \ldots, p_n, q_n)
\end{gather*}
having the coef\/f\/icients of $p_r$ and $q_r$ as column vectors.
One then def\/ines the quantity
\begin{gather*}
P = \prod\limits_r \det(\mathbf{u}^{-}_{r},\mathbf{u}^{+}_r) \prod\limits_{r<s}
\left(\det(\mathbf{u}^{-+}_{rs},\mathbf{u}^{+-}_{rs}) \det(\mathbf{u}^{--}_{rs},\mathbf{u}^{++}_{rs}) \right)^2.
\end{gather*}
Then the determinant function
\begin{gather*}
D(x_1,\ldots,x_n) = \det(M)/P
\end{gather*}
is independent of the choices of lifts, and is thus well def\/ined.
Similar to the Atiyah--Sutclif\/fe determinant,~$D$ is actually invariant under the action of the Weyl group $W_n$ on
$C_n$, and is also invariant under scaling, and rotations in $\mathbb{R}^3$.
However, unlike the Atiyah--Sutclif\/fe determinant, it is not invariant under translations, because the origin in
$\mathbb{R}^3$ plays here a~special role, and it is always real-valued, because it is the determinant of a~$2n$ by $2n$
complex matrix, which represents an~$n$ by~$n$ quaternionic matrix, and thus is always real (indeed, the complex
conjugate of such a~$2n$ by $2n$ complex matrix can be shown to be in the same conjugacy class as the complex matrix
itself, so they must have equal determinants).

\section[The case $n=2$]{The case $\boldsymbol{n=2}$}

We consider here the case $n=2$.
We have two points $x_1, x_2 \in \mathbb{R}^3$ such that $x_1 \neq x_2$ and $x_1 \neq -x_2$.
Using a~rotation in $\mathbb{R}^3$, we can assume that they both lie on the $xy$-plane.
We think of the $xy$-plane as the complex plane.
Using a~rotation in the $xy$-plane and scaling, we can further assume that $x_1=1$ and we then let~$z$ be the complex
number representing $x_2$ in the $xy$-plane.
Thus $z\neq 1$ and $z \neq -1$.
We let
\begin{gather*}
A = \frac{z-1}{|z-1|},
\qquad
B = -\frac{z+1}{|z+1|},
\qquad
g = -\frac{z}{|z|}.
\end{gather*}
We then have
\begin{gather*}
-64 D = \left|
\begin{matrix} AB & 1 & ABg & 1
\\
AB-A-B & \bar{A}+\bar{B}-1 & -Ag+Bg-AB & -\bar{A}+\bar{B}+\bar{g}
\\
1-A-B & \bar{A}\bar{B}-\bar{A}-\bar{B} & A-B-g & -\bar{A}\bar{g}+\bar{B}\bar{g}-\bar{A}\bar{B}
\\
1 & -\bar{A}\bar{B} & 1 & -\bar{A}\bar{B}\bar{g}
\end{matrix}
\right|.
\end{gather*}
We then multiply the second column by $-AB$ and add it to the f\/irst column, and we multiply the second column by $-ABg$
and add it to the third column, and f\/inally subtract the second column from the fourth one, and get, after expanding the
determinant along the f\/irst row:
\begin{gather*}
64 D = \left|
\begin{matrix} 2(AB-A-B) & -2Ag+ABg-AB & -2\bar{A}+1+\bar{g}
\\
0 & -2g+Bg-B+Ag+A & -2\bar{A}\bar{B}+\bar{A}(1-\bar{g})+\bar{B}(1+\bar{g})
\\
2 & 1+g & \bar{A}\bar{B}(1-\bar{g})
\end{matrix}
\right|.
\end{gather*}
Taking a~$2$ out from the f\/irst column, and using elementary column operations using the f\/irst column in order to make
the entries in the $(3,2)$ and $(3,3)$ positions vanish, we get
\begin{gather*}
32 D  = \left|
\begin{matrix}
AB-A-B & -2AB+A(1-g)+B(1+g) & 2\bar{g}-\bar{A}(1+\bar{g})+\bar{B}(1-\bar{g})
\\
0 & A(1+g)-B(1-g)-2g & -2\bar{A}\bar{B}+\bar{A}(1-\bar{g})+\bar{B}(1+\bar{g})
\\
1 & 0 & 0
\end{matrix}
\right|
\\
\phantom{32 D}{}
 = |-2AB+A(1-g)+B(1+g)|^2+|A(1+g)-B(1-g)-2g|^2
\\
\phantom{32 D}{}
 = 8 + 2|1-g|^2 + 2|1+g|^2-2B(1-\bar{g})-2 \bar{B}(1-g) -2A(1+\bar{g}) -2\bar{A}(1+g)+\cdots
\\
\phantom{32 D=}{}
 +2(A\bar{B}-\bar{A}B)(\bar{g}-g) -2A\bar{g}(1+g)-2\bar{A}g(1+\bar{g})+2B\bar{g}(1-g)+2\bar{B}g(1-\bar{g}).
\end{gather*}
Using
\begin{gather*}
|1+g|^2 + |1-g|^2 = 4
\end{gather*}
we get
\begin{gather*}
16 D = 8-2A(1+\bar{g})-2B(1-\bar{g})+A(1+\bar{g})\bar{B}(1-g)-2\bar{A}(1+g)-2\bar{B}(1-g)
\\
\phantom{16 D =}{}
+\bar{A}(1+g)B(1-\bar{g}).
\end{gather*}
If we let
\begin{gather*}
(w_1,w_2)  = \frac{1}{2}(w_1 \bar{w}_2 + w_2 \bar{w}_1),
\qquad
\det(w_1,w_2)  = \frac{i}{2}(w_1\bar{w}_2-w_2\bar{w}_1),
\end{gather*}
we can then write
\begin{gather*}
4D = 2-(A,1+g)-(B,1-g)-\Im(g)\det(A,B),
\end{gather*}
where $\Im(g)$ denotes the imaginary part of~$g$.
Therefore, using the def\/initions of~$A$,~$B$ and~$g$ in terms of~$z$, we get
\begin{gather*}
4D = 2+\left(\frac{z-1}{|z-1|},\frac{z}{|z|}-1\right)+\left(\frac{z+1}{|z+1|},\frac{z}{|z|}+1\right)+2\frac{(\Im{z})^2}{|z||z-1||z+1|}.
\end{gather*}
Writing $z=re^{i\theta}$, and after simplif\/ication, we get
\begin{gather*}
4D = 2 + \frac{(1+r)(1-\cos(\theta))}{|z-1|} + \frac{(1+r)(1+\cos(\theta))}{|z+1|} +
\frac{2r(1-\cos(\theta))(1+\cos(\theta))}{|z-1||z+1|}.
\end{gather*}
Using $1+r \geq |z+1|$ and $1+r \geq |z-1|$, and that $1+\cos(\theta)$ and $1-\cos(\theta)$ are both nonnegative,
\begin{gather*}
4D \geq 4 + \frac{2r(1-\cos(\theta))(1+\cos(\theta))}{|z-1||z+1|}.
\end{gather*}
Thus
\begin{gather*}
D \geq 1 + \frac{r\sin^2(\theta)}{2|z-1||z+1|}.
\end{gather*}
This proves the inequality $D \geq 1$, which in turns implies the linear independence conjecture, for $n=2$.
Moreover, it is not too dif\/f\/icult to see that equality $D=1$ occurs if and only if $\sin(\theta) = 0$, or, in other
words, if the two points $x_1$ and $x_2$ lie on the same line through the origin.

\section{A stronger conjecture}

Similar to Conjecture~2 in~\cite{Atiyah-Sutcliffe2002}, we make the following conjecture

\begin{conjecture}
For any $n\geq 2$ and for any $\mathbf{x} \in C_n$, we have $D(\mathbf{x}) \geq 1$.
\end{conjecture}

The author did some numerical testing for this conjecture for $n \leq 10$, by computing $D(x)$ for a~number of
pseudo-randomly generated conf\/igurations of points, and found that in all these cases, the conjecture was verif\/ied, thus
gathering some numerical evidence for the above conjecture.

\subsection*{Acknowledgements}
The author would like to thank Sir Michael Atiyah for kindly replying to his emails, and would like to thank the
anonymous referees for all their suggestions, which ended up making the article much more readable.

\pdfbookmark[1]{References}{ref}
\LastPageEnding

\end{document}